\documentclass[12pt,reqno]{amsart}

\DeclareMathOperator{\Gal}{Gal}%

\DeclareMathOperator{\Br}{Br}

\usepackage{amssymb}
\usepackage{hyperref}
\usepackage[arrow]{xypic}

\begin{document}

\title[Galois module structure of Milnor $K$-theory]{Galois
module structure of Milnor $K$-theory mod $p^s$ in characteristic
$p$}

\author[J\'{a}n Min\'{a}\v{c}]{J\'an Min\'a\v{c}$^{*}$}
\address{Department of Mathematics, Middlesex College, \ University
of Western Ontario, London, Ontario \ N6A 5B7 \ CANADA}
\thanks{$^*$Research supported in part by NSERC grant R0370A01, and by a Distinguished Research Professorship at the University of Western Ontario.}
\email{minac@uwo.ca}

\author[Andrew Schultz]{Andrew Schultz}
\address{Department of Mathematics, University of Illinois
at Urbana-Cham\-paign, 273 Altgeld Hall, MC-382, 1409 W.~Green
Street, Urbana, IL \ 61801 \ USA}
\email{acs@math.uiuc.edu}

\author[John Swallow]{John Swallow$^\dag$}
\address{Department of Mathematics, Davidson College, Box 7046, Davidson, North Carolina \ 28035-7046 \ USA}
\thanks{$^\dag$Research supported in part by NSF grant DMS-0600122.}
\email{joswallow@davidson.edu} 

\begin{abstract}
    Let $E$ be a cyclic extension of $p$th-power degree of a field $F$ of characteristic $p$.  For all $m$, $s\in \mathbb{N}$, we determine $K_mE/p^sK_mE$ as a $(\mathbb{Z}/p^s\mathbb{Z}) [\Gal(E/F)]$-module.  We also provide examples of extensions for which all of the possible nonzero summands in the decomposition are indeed nonzero.
\end{abstract}


\date{September 5, 2007}

\maketitle

\newtheorem*{theorem*}{Theorem}
\newtheorem{corollary}{Corollary}
\newtheorem{proposition}{Proposition}
\newtheorem{lemma}{Lemma}

\theoremstyle{definition}
\newtheorem*{definition*}{Definition}
\newtheorem*{remarks*}{Remarks}

\newcommand{\B}{B}
\newcommand{\Bb}{\mathcal{B}}
\newcommand{\F}{\mathbb{F}}
\newcommand{\Fp}{\F_p}
\newcommand{\N}{\mathbb{N}}
\newcommand{\R}{\mathbb{R}}
\newcommand{\Z}{\mathbb{Z}}
\newcommand{\comment}[1]{}

\parskip=12pt plus 2pt minus 2pt


Let $F$ be a field of characteristic $p$.  Let $K_mF$ denote the $m$th Milnor $K$-group of $F$ and $k_mF=K_mF/pK_mF$. (See, for instance, \cite{Mi} and \cite[Chapter~14]{Ma}.) If $E/F$ is a Galois extension of fields, let $G=\Gal(E/F)$ denote the associated Galois group.  In \cite{BLMS} the structure of $k_mE$ as an $\Fp G$-module was determined when $G$ is cyclic of $p$th-power order.  In this paper we determine the Galois module structure of $K_mE$ modulo $p^s$ for $s \in \N$ and these same $G$. We also provide examples of extensions for which the possible free summands in the decomposition are all nonzero. These examples together with the results in \cite{BLMS} show that the dimensions over $\Fp$ of indecomposable $\Fp[\Gal (E/F)]$-modules occuring as direct summands of $k_mE$ are all powers of $p$ and that all dimensions $p^i$, $i=0, 1, \dots,n$, indeed occur in suitable examples.

Recall the theorem of Bloch-Kato and Gabber (see \cite{BK}): the sequence
\begin{equation*}
    0 \to k_m F \to \Omega_F^m \stackrel{\mathfrak{P}}
    {\to} \Omega_F^m/d \; \Omega_F^{m-1}
\end{equation*}
is exact, where $\Omega_F^m$ is the $m$th graded component of the exterior algebra on K\"ahler differentials and $\mathfrak{P}$ is the Artin-Schreier operator. In \cite[\S 6]{I}, Izhboldin succeeded in providing an analogue of this important interpretation of $k_m F$, as follows: for $s \in \N$, the sequence
\begin{equation*}
    0 \to K_m F/p^s K_m F \stackrel{\delta} {\to} Q^m (F,s)
    \stackrel{\mathfrak{P}}{\to} Q^m (F,s)
\end{equation*}
is exact, where $Q^m (F,s)$ is the ``Milnor-Witt'' group of $F$, defined together with maps $\delta$ and $\mathfrak{P}$ in \cite[\S 6]{I}. These objects play an important role in the arithmetic of fields, higher class field theory, and Milnor $K$-theory of fields of characteristic $p$.  (See, for instance, \cite{FK} and the references therein.)

Since the classification problem of $(\Z/p^s\Z)G$-modules for cyclic $G$ is nontrivial and has not been completely solved---see, for instance, \cite{T} for results and references---it is a pleasant surprise that the $(\Z/p^s\Z )G$-modules $K_m E/p^s K_m E$ have a simple description. The main ingredients we use to obtain this description are the lack of $p$-torsion in $K_m E$, due to Izhboldin \cite{I}, together with the result \cite{BLMS} for the case $s=1$---which also depends on Izhboldin's result.

Suppose that $E/F$ is cyclic of degree $p^n$, and for $i=0, 1, \dots, n$, let $E_i/F$ be the subextension of degree $p^i$ of $E/F$ and $G_i:=\Gal(E_i/F)$.  Set $R_s:=\Z/p^s\Z$ and let $\Z_p$ be the ring of $p$-adic integers. We write $\iota_{F,E}\colon K_mF\to K_mE$ and $N_{E/F}\colon K_mE\to K_mF$ for the natural inclusion and norm maps, and we use the same notation for the induced maps between $K_mF/p^sK_mF$ and $K_mE/p^sK_mE$.

\begin{theorem*}
    There exists an isomorphism of $R_sG$-modules
    \begin{equation*}
        K_mE/p^sK_mE \simeq \oplus_{i=0}^n Y_i
    \end{equation*}
    where
    \begin{itemize}
        \item $Y_n$ is a free $R_sG$-module of rank $\dim_{\Fp}
        N_{E/F}k_mE$,
        \item $Y_i$, $0<i<n$, is a free $R_sG_i$-module of rank
        \begin{equation*}
            \dim_{\Fp} N_{E_i/F} k_mE_i/N_{E_{i+1}/F}k_mE_{i+1},
        \end{equation*}
        \item $Y_0$ is a free $R_s$-module of rank
        \begin{equation*}
            \dim_{\Fp} k_mF/N_{E_1/F} k_mE_1 \mbox{ with }
        Y_0^G=Y_0,
        \end{equation*}
        and
        \item for each $0\le i\le n$, $Y_i\subseteq
        \iota_{E_i,E}(K_mE_i/p^sK_mE_i)$.
    \end{itemize}
    Moreover,
    \begin{equation*}
        \widehat{K_m E} := \varprojlim_s K_m E/p^s K_m E
        \simeq \oplus_{i=0}^n \; \hat Y_i,
    \end{equation*}
    where each $\hat Y_i$ is a free $\Z_p G_i$-module of the same
    rank as $Y_i$.
\end{theorem*}

\section{Proof of the Theorem}

    We prove the result by induction on $s$.  The case $s=1$ is \cite[Theorem~2]{BLMS}.  Assume therefore that $s>1$ and the result holds for $s-1$:
    \begin{equation*}
        K_mE/p^{s-1}K_mE = \oplus \tilde Y_i,
    \end{equation*}
    with each $\tilde Y_i$ a free $R_{s-1}G_i$-module $\tilde Y_i$ in the image of $\iota_{E_i,E}$, $0\le i\le n$.

    For each $i$ with $0\le i\le n$, let $\B_{s-1,i} \subset i_{E_i,E} K_mE_i/p^{s-1} K_mE_i$ be an $R_{s-1}G_i$-base for the free $R_{s-1}G_i$-module $\tilde Y_i$. By induction the cardinality of $\B_{s-1,i}$ is
    \begin{equation*}
        \vert \B_{s-1,i}\vert =
        \begin{cases}
        \dim_{\Fp} N_{E_i/F}k_mE_i/N_{E_{i+1}/F} k_mE_{i+1},
        & i<n \\
        \dim_{\Fp} N_{E/F}k_mE, & i=n.
        \end{cases}
    \end{equation*}
    Since the $\tilde Y_i$ are independent, the set
    \begin{equation*}
        \B_{s-1} := \cup_{0\le i\le n} \B_{s-1,i} \subseteq
        K_mE/p^{s-1}K_mE
    \end{equation*}
    is $R_{s-1}G$-independent.

    For each $i$, let $\Bb_i\subseteq \iota_{E_i,E}(K_mE_i)$ be a set of representatives for the elements of $\B_{s-1,i}$, and let $\B_{s,i}\subseteq K_mE/p^sK_mE$ be chosen to make the following first diagram commutative. The second diagram merely recalls where our $\Bb_{i},B_{s,i}$ and $B_{s-1,i}$ are located.
    \begin{equation*}
        \xymatrix{
        \Bb_{i} \ar[r]^-{\bmod p^{s}} \ar[d]_-{\bmod p^{s-1}} & \B_{s,i}
        \ar[dl]^{\bmod p^{s-1}} & & K_mE \ar[r] \ar[d] &
        K_mE/p^{s}K_mE \ar[dl] \\ \B_{s-1,i} & & & K_mE/p^{s-1}K_mE
        }
    \end{equation*}
    Hence for each $i$ we have bijections
    \begin{equation*}
        \Bb_i \leftrightarrow \B_{s,i}\leftrightarrow
        \B_{s-1,i}
    \end{equation*}
    and $\vert \B_{s,i}\vert = \vert \B_{s-1,i}\vert$.

    First we observe that every nonzero ideal $V$ of $R_sG_i$
    contains $p^{s-1}$ $(\tau-1)^{p^i}$, where $\tau$ is any
    fixed generator of $G_i$. Indeed consider $0 \ne \beta
    \in B$. By multiplying by an appropriate power of $p$, we may
    assume $0 \neq \beta \in p^{s-1}R_sG_i$.  Let us write
      \begin{equation*}
        \beta = \sum_{j=k}^{p^i-1} c_j (\tau-1)^j,
      \end{equation*}
    where each $c_j \in p^{s-1}R_s$, $j = k, \dots, p^i-1$, and $c_k
    \not\in p^sR_s=\{0\}$, say $c_k = p^{s-1} \tilde c_k$ with
    $\tilde c_k \not\in pR_s$. Using the fact that $p^{s-1}
    (\tau-1)^{p^i} = 0$ in $R_sG_i$ we see that we can multiply
    $\beta$ by $(\tau-1)^{p^i-k-1}$ to obtain
      \begin{equation*}
        0 \ne \tilde c_k p^{s-1} (\tau-1)^{p^i-1} \in V.
      \end{equation*}
    Since $\tilde c_k \in U(R_s)$ we see that $p^{s-1}
    (\tau-1)^{p^i-1} \in V$ as asserted.

    Set $Y_i$ to be the $R_sG$-submodule of $K_mE/p^sK_mE$ generated by $\B_{s,i}$.  It is clear that $Y_i\subseteq \iota_{E_i,E}(K_mE_i/p^sK_mE_i)$ and hence $Y_i$ is an $R_sG_i$-module.

    Each element $b\in \B_{s,i}$ generates in $K_mE/p^sK_mE$ a free $R_sG_i$-module $M_b$, as follows.  Suppose that $M_b$ is not a free $R_sG_i$-module. Then the annihilator of $b$ in $R_sG_i$ is a nonzero ideal of $R_sG_i$. Let $\hat b\in \Bb_i$ and $\tilde b\in \B_{s-1,i}$ correspond to $b$ under the bijection above. Let also $\sigma$ be a generator of $G$ and $\bar\sigma$ its image in $G_i$. Since every nonzero ideal of $R_sG_i$ contains $p^{s-1}(\bar\sigma-1)^{p^i-1}$, for some $\hat c\in K_mE$ we have
    \begin{equation*}
        p^{s-1}(\bar\sigma-1)^{p^i-1}\hat b = p^s\hat c.
    \end{equation*}

    Since $K_mE$ has no $p$-torsion \cite[Theorem A]{I}, we obtain
    \begin{equation*}
        p^{s-2}(\bar\sigma-1)^{p^i-1}\hat b = p^{s-1}\hat c.
    \end{equation*}
    Then in $K_mE/p^{s-1}K_mE$
    \begin{equation*}
        p^{s-2}(\bar\sigma-1)^{p^i-1}\tilde b = 0,
    \end{equation*}
    contradicting the fact that $\tilde b$ lies in the $R_{s-1}G_i$-base $\B_{s-1,i}$ for $\tilde Y_i$. (Alternatively we could use \cite[Theorem~5.1]{T} to show that $M_b$ is a free $R_s G_i$-module.)

    Now set $\B_s:=\cup_{0\le i\le n} \B_{s,i} \subseteq K_mE/p^sK_mE$.  Suppose we have a relation
    \begin{equation*}
        \sum r_\alpha b_\alpha = 0, \quad\quad b_\alpha\in \B_s,\ \
        r_\alpha\in R_sG_\alpha,
    \end{equation*}
    such that $G_\alpha = G_i$ for a suitable $i$, $0\le i\le n$, with $b_\alpha \in B_{s,i}$. Let $\tilde b_\alpha \in \B_{s-1}$ correspond to $b_\alpha$ under the natural bijection, and similarly let $\tilde r_\alpha \in R_{s-1} G_\alpha$ be the image of $r_\alpha \in R_s G_\alpha$.

    Working mod $p^{s-1}$ we have $\sum \tilde r_\alpha \tilde b_\alpha = 0$. Since each $\tilde b_\alpha$ lies in the $R_{s-1}G$-independent set $\B_{s-1}$, we deduce that $r_\alpha\in p^{s-1}R_sG_\alpha$ for each $\alpha$. Write $r_\alpha = ps_\alpha$ for elements $s_\alpha \in R_s G_\alpha$. We rewrite the original relation as
    \begin{equation*}
        \sum ps_\alpha b_\alpha = 0, \quad\quad b_\alpha\in \B_s,
        \ \ s_\alpha\in R_s G_\alpha.
    \end{equation*}
    Just as before we divide by $p$ to obtain in $K_mE/p^{s-1}K_mE$
    \begin{equation*}
        \sum \tilde s_\alpha \tilde b_\alpha = 0.
    \end{equation*}
    Again since each $\tilde b_\alpha\in \B_{s-1}$, we deduce that $s_\alpha\in p^{s-1}R_s G_\alpha$. But then $r_\alpha = ps_\alpha = 0\in R_s G_\alpha$ for each $r_\alpha$, as desired.

    Hence for each $i$ in $0\le i\le n$ we have that $Y_i$ is a direct sum of free $R_sG_i$-modules $M_b$ for $b\in \B_{s,i}$, and moreover that $\sum Y_i = \oplus Y_i$.  By Nakayama's Lemma, since $\Bb$ generates $K_m/p^{s-1}K_mE$ it also generates $K_m/p^sK_mE$, and hence $\oplus Y_i = K_mE/p^sK_mE$.  (More explicitly, choose $\alpha \in K_mE/p^sK_mE$ and $\hat \alpha \in K_mE$ a lift of $\alpha$. Since $B_{s-1}$ spans $K_mE/p^{s-1}K_mE$ we have
    \begin{equation*}
          \hat \alpha = \sum_{b \in \Bb} f_b b + p^{s-1}\hat \gamma
    \end{equation*}
    for some $\hat \gamma \in K_mE$, where each $f_b \in R_sG$ and all but finitely many $f_b =0$.  We also have $\hat \gamma = \sum g_b b+p^{s-1}\hat \delta$ for some $\hat \delta \in K_mE$, where again each $g_b \in R_sG$ and all but finitely many $g_b =0$. Therefore
    \begin{equation*}
    \alpha = \sum_{b \in \Bb} \left(f_b+p^{s-1}g_b\right)b
    \end{equation*}
    as elements of $K_mE/p^sK_mE$.)

    The last statement concerning the equality $\widehat{K_m E} = \oplus_{i=0}^n \; \hat Y_i$ is obtained by passing to projective limits.\qed

\section{Examples of $E/F$ with $Y_i\neq \{0\}$ for all $i$}

Let $p$ be an arbitrary prime number, and let $q$ be an arbitrary prime number or $0$.  We show that for each $n,m \in \N$ there exists a cyclic field extension $E/F$ of degree $p^n$ and characteristic $q$ such that for each $i$, $0\leq i < n$,
\begin{equation*}
    \dim_{\Fp} N_{E_i/F} k_m E_i/N_{E_{i+1}/F} k_m E_{i+1} \neq 0
\end{equation*}
and
\begin{equation*}
    \dim_{\Fp} N_{E/F} k_m E \neq 0.
\end{equation*}
Recall that we index $E_i$ such that $F \subset E_i \subset E$
and $[E_i \colon F] = p^i$.

We are interested in the images of the quotient groups in $k_mE$. Because in the case $p=q$ the natural homomorphism $k_mF \to k_m E$ is injective, for the case of our main application we can work in $k_mF$.

\subsection{The case $m=1$}\

Fix $p$, $q$, and $n$ as above and set $m=1$.  We construct a field extension $E/F$ as above together with elements
\begin{align*}
    x_i&\in N_{E_i/F}(E_i^\times)\setminus N_{E_{i+1}/F}(E_{i+1}^\times)F^{\times p},\quad 0\le i< n, \\
    x_n&\in N_{E/F}(E^\times)\setminus F^{\times p}.
\end{align*}

Let $B$ be a field of characteristic $q$, and let $A/B$ be a cyclic extension of degree $p^n$.  Index the subfields $A_i$ of $A/B$ such that $[A_i:B] = p^i$, and denote by $\iota_{B,A_i}: K_1B \hookrightarrow K_1A_i$ the natural inclusion. Let $\sigma$ be a generator of $\Gal(A/B)$ and set $\sigma_i = \sigma \vert_{A_i}$, the restricted map. Finally, assume that there exist elements $x_0, x_1,\dots,x_n \in B^\times$ such that the following condition (\verb+*+) holds:
\begin{equation*}\label{eqn:condition}
\begin{split}
    [\iota_{B,A_j}(x_{j})]^{p^{n-j-1}} &\notin \langle
    [\iota_{B,A_j}(x_1)]^{p^{n-1}},
    [\iota_{B,A_j}(x_2)]^{p^{n-2}},\dots,
    [\iota_{B,A_j}(x_{n})] \rangle \\  &\subseteq
    A^\times_{j}/N_{A/A_{j}} (A^\times), \quad 0\le j<n,\\ x_n
    &\notin B^{\times p},
\end{split}
\end{equation*}
where $[x]$ denotes the class of $x$ and $\langle S\rangle$ the subgroup generated by a set $S$ in the named factor group. At the end of this section we shall create an example where condition (\verb+*+) holds.

Now consider cyclic algebras
\begin{equation*}
    \mathcal{A}_{j} = (A/B,\sigma, x_{j}^{p^{n-j}}),
    \quad 1\le j\le n.
\end{equation*}
Observe that
\begin{equation*}
    [\mathcal{A}_j] = [(A_{j}/B, \sigma_{j},
    x_{j})] \in \Br(B), \quad 1\le j\le n
\end{equation*}
(\cite[Chapter~15, Corollary~b]{P}), where $\Br(B)$ denotes the Brauer group of $B$.  Let $F$ be the function field of the product of Brauer-Severi varieties of $\mathcal{A}_1,\dots,\mathcal{A}_{n}$. (See \cite[page~735]{SV}; see also \cite[Chapter~3]{J} for basic properties of Brauer-Severi varieties.)

Let $E$ be the compositum $A\cdot F$ of the fields $F$ and $A$. Since $F$ is a regular extension of $B$, we see that $E/F$ is a cyclic extension of degree $p^n$. We denote again as $\sigma$ the generator of $\Gal(E/F)$ which restricts to $\sigma \in \Gal(A/B)$, and we write $E_k = A_k\cdot F$ for $k = 0, 1, \dots,n$. Now $[\mathcal{A}_j \otimes_B F] = 0 \in \Br(F)$, $j=1, \dots,n$, because $F$ splits each $\mathcal{A}_j$. Hence
\begin{equation*}
        0 = [(E/F, \sigma, x_j^{p^{n-j}})] =
        [(E_{j}/F, \sigma_{j}, x_j)],
\end{equation*}
and so $x_j \in N_{E_{j}/F} (E^\times_{j})$ as desired (see
\cite[Chapter~15, page~278]{P}).

However, we claim that
\begin{align*}
    x_j &\notin (N_{E_{j+1}/F} (E^\times_{j+1})) F^{\times p}, \quad
    0\le j<n, \\ x_n &\notin F^{\times p}.
\end{align*}
Since $x_n\notin B^{\times p}$ by hypothesis and $F/B$ is a regular extension, we have $x_n\notin F^{\times p}$.  Assume then that $0\le j<n$ and, contrary to our statement,
\begin{equation*}
    x_j \in (N_{E_{j+1}/F}( E^\times_{j+1}))F^{\times p}.
\end{equation*}
Then we have $x_j f^p \in N_{E_{j+1}/F} (E^\times_ {j+1})$ for some $f \in F^\times$.  Hence
\begin{equation*}
    [(E_{j+1}/F, \sigma_{j+1}, x_j f^p)] = 0 \in \Br (F)
\end{equation*}
and so
\begin{align*}
    [(E_{j+1}/F, \sigma_{j+1}, x_j)]
    &=-[(E_{j+1}/F, \sigma_{j+1}, f^p)] \\
    &=-[(E_{j}/F, \sigma_{j}, f)].
\end{align*}
(In the case $j=0$, we use $(E_0/F,\sigma_0,f)$ to denote the zero element in $\Br(F)$.) Consequently $(E_{j+1}/F, \sigma_{j+1}, x_j)$ is split by $E_{j}$. (See \cite[Chapter~15, Proposition~b]{P}.)

But then
\begin{equation*}
    [(E_{j+1}/E_{j}, \sigma_{j+1}^ {p^{j}}, \iota_{F,E_j}(x_j))] =
    0 \in \Br(E_{j}).
\end{equation*}
(See \cite[page~74]{D}.) Hence $[(E/E_{j}, \sigma^{p^{j}}, \iota_{F,E_j}(x_j^{p^{n-j-1}}))] = 0 \in \Br(E_{j})$. But $E_{j}=A_j\cdot F$ is the function field of the product of the Brauer-Severi varieties of $\mathcal{A}_k \otimes_B A_{j}$ for $k = 1,\dots,n$. Therefore
\begin{equation*}
    [(A/A_{j}, \sigma^{p^{j}}, \iota_{B,A_j}(x_j^ {p^{n-j-1}}))] \in
    \langle [\mathcal{A}_k \otimes_B A_{j}], k=1, \dots, n \rangle
    \subseteq \Br (A_{j}).
\end{equation*}
Consequently
\begin{align*}
    [\iota_{B,A_j}(x_j^{p^{n-j-1}})] &\in \langle
    [\iota_{B,A_j}(x_1)]^{p^{n-1}},\dots,[\iota_{B,A_j}(x_n)]\rangle
    \\ & \in A_j^\times/N_{A/A_{j}}(A^\times),
\end{align*}
a contradiction to condition (\verb+*+).

Thus we have shown that a required extension $E/F$ exists with elements $x_0, x_1, \dots,x_n$, provided that we can produce a field extension $A/B$ and elements $x_0, x_1, \dots, x_n \in B^\times$ such that condition (\verb+*+) is valid.  Now we show that such an extension and elements exist.

Let $B := C (x_0, x_1, \dots, x_n)$, where $C$ is a field with characteristic $q$ and $x_0, x_1, \dots,x_n$ are algebraically independent elements over $C$. Assume also that there exists a cyclic extension $D/C$ of degree $p^n$ with Galois group $G = \langle \sigma \rangle$. Finally, let $A := D (x_0, x_1, \dots, x_n)$. Thus $A/B$ is a cyclic extension of degree $p^n$.

We claim that condition (\verb+*+) holds. Clearly $x_n\notin
B^{\times p}$. Contrary to our claim assume that
\begin{equation*}
    x_j^{p^{n-j-1}} = x_1^{c_1p^{n-1}}\cdots x_{n}^{c_{n}}
    N_{A/A_{j}}(\gamma)
\end{equation*}
where $c_i \in \Z$, $\gamma \in A^\times$, and $0\le j<n$. Write
\begin{equation*}
    \gamma = u \cdot \frac{p_1 \cdots p_t}
    {q_1 \cdots q_w},
\end{equation*}
where $\{ p_1, \dots, p_t \}$ and $\{ q_1, \dots, q_w \}$ are disjoint sets of primes in the ring $D [x_0, x_1, \dots,x_n]$ and $u\in U(D
[x_0, x_1, \dots, x_n])=D^\times$.

We write $H_j = \Gal (A/A_{j})$ for $j = 0,1,\dots, n-1$.  Then we
have
\begin{equation*}
    x_j^{p^{n-j-1}} \prod_{h \in H_j}\prod_{i=1}^w
    h (q_i) = x_1^{c_1p^{n-1}} \cdots
    x_{n}^{c_{n}} N_{A/A_j} (u)
    \prod_{h \in H_j} \prod_{i=1}^t h (p_i).
\end{equation*}

In order to show that this equation is impossible, consider the discrete valuation $v_j$ on $A$ such that $v_j (x_j) = 1$, $v_j(x_k) = 0$ if $j \ne k$, and $v_j(d)=0$ for $d \in D$. Consider first the case when $x_j \nmid q_l$, $l = 1,\dots,w$. Then the value of $v_j$ on the left-hand side is $p^{n-j-1}$, while the value of $v_j$ on the right-hand side is at least $p^{n-j}$. Indeed if $c_j \ne 0$ this is true as $x_j^{p^{n-j}}$ divides the right-hand side of our equation. If $c_j = 0$ or $j=0$ (in which case $c_0$ is not defined) then since $x_j$ divides the left-hand side we see that there exists $p_l$, $l \in \{ 1,\dots,t \}$, such that $x_j \mid p_l$. Because $\vert H_j \vert = p^{n-j}$ we see again that the value of $v_j$ on the right-hand side is again at least $p^{n-j}$. Thus $x_j \mid q_l$ for some $l = 1,\dots,w$, and so $x_j \nmid p_l$ for any $l=1,\dots,t$. Then the value of $v_j$ on the left-hand side of our equation is $p^{n-j-1} + p^{n-j}$ while the value of $v_j$ on the right-hand side is at most $p^{n-j}$. Thus we see that the equation above is impossible. Hence condition (\verb+*+) is valid and we have established the desired example of $E/F$ in the case $m = 1$.

\subsection{The case $m>1$}\

Fix $m$, $n$, and $q$ as above and let $L/K$ be a field extension satisfying the case $m = 1$ with the same $n$ and $q$.  Let $x_0, x_1,\dots, x_n \in K^\times$ such that
\begin{align*}
    x_i&\in N_{L_i/K}(L_i^\times)\setminus
    N_{L_{i+1}/K}(L_{i+1}^\times)K^{\times p},\quad 0\le i< n, \\
    x_n&\in N_{L/K}(L^\times)\setminus K^{\times p}.
\end{align*}
Consider the field of the iterated power series $F := K ((y_1)) \cdots ((y_{m-1}))$. Then $E:=L\cdot F$ is a cyclic extension of degree $p^n$ over $F$. For each $j \in \{ 0,1,\dots,n \}$ consider the element
\begin{equation*}
    \alpha_j = \{ x_j, y_1,\dots,y_{m-1} \} \in
    k_m F.
\end{equation*}
(If $m = 1$ then $\alpha_j = \{ x_j \}$.) By our hypothesis and the projection formula for the norm map in $K$-theory, we have
\begin{equation*}
    \alpha_j \in N_{E_{j}/F} k_m E_{j}.
\end{equation*}

Now for $0\le j<n$ we shall prove by induction on $m \in \N$ that
    \begin{equation*}
        \alpha_j \notin N_{E_{j+1}/F}
        k_m E_{j+1}.
\end{equation*}
If $m=1$ then our statement is true by the choice of the field extension $L/K$ and the elements $x_j$. Assume then that $m > 1$ and that our statement is true for $m - 1$.

Consider the complete discrete valuation $v$ on $F$ with uniformizer $y_{m-1}$ and residue field $F_v = K ((y_1)) \cdots ((y_{m-2}))$ if $m>2$ and $F_v = K$ if $m=2$.

For the sake of simplicity we denote by $E'$ the field $E_{j+1}$, and denote the unique extension of $v$ on $E'$ again by $v$. Since we are considering an unramified extension we assume that both valuations are normalized. Let $\partial: k_m F \to k_{m-1} F_v$ and $k_m E' \to k_{m-1} E'_v$ be the homomorphisms induced by residue maps in Milnor $K$-theory. Then applying \cite[Lemma~3]{K} we see that the following diagram is commutative:
    \begin{equation*}
        \xymatrix{k_m E' \ar[r]^-\partial
        \ar[d]_{N_{E'/F}} & k_{m-1} E'_v \ar[d]^
        {N_{E'_v/F_v}} \\ k_m F \ar[r]^-\partial &
        k_{m-1} F_v.}
    \end{equation*}
If $\alpha_j \in N_{E'/F} k_m E'$, then $\partial
\alpha_j \in N_{E'_v/F_v} k_{m-1} E'_v$.

But $\partial \alpha_j = \{ x_j, y_1,\dots,y_{m-2} \}$ if $m > 2$  and $\partial \alpha_j = \{x_j\}$ if $m = 2$, a contradiction in either case.  Therefore we have constructed a field extension $E/F$ with the desired properties.

\begin{remarks*}
    In \cite{MSS1} we determined the $\Fp G$-module structure of $k_1 E$ for all cyclic extensions $E/F$ of degree $p^n$, where $G$ is the Galois group.  In particular, the decomposition does not depend upon the characteristic of the base field. The ranks of the free $\Fp G_i$-summands appearing in that decomposition are again determined by the image of $N_{E_i/F} (E_i^\times)/N_{E_{i+1}/F} (E^\times_{i+1})$ in $E^\times/ E^{\times p}$.

    When no primitive $p$th root of unity lies in $F$, we have that $F^\times/F^{\times p}$ embeds in $E^\times/E^{\times p}$. Therefore the construction given above for field extensions $E/F$ applies in this case, and these free $\Fp G_i$-summands do indeed occur for any characteristic. When a primitive $p$th root of unity is in $F^\times$, the kernel of the homomorphism $F^\times/ F^{\times p} \to E^\times/E^{\times p}$ is generated by a class $[a] \in F^\times/F^{\times p}$, where $E_1 = F (\root{p}\of{a})$. In this case it is enough to additionally require that $a \in N_{E/F} (E^\times)$ when $p=2$ and $n=1$ (since the condition is automatic otherwise), and the construction above of field extensions $E/F$ applies again.

    If $F$ contains a primitive $p$th-root of unity and $m=1$, the decomposition contains at most one other indecomposable module, a cyclic $\Fp G$-module of dimension $p^k + 1$ for $k \in \{ -\infty, 0, 1, \dots, n-1 \}$ (where we set $p^{-\infty}=0$). In \cite{MSS2} we showed that all of these modules are realizable as well.
\end{remarks*}

\section*{Acknowledgments}

Andrew Schultz would like to thank Ravi Vakil for his encouragement and direction in this and in all other projects. We would also like to thank Ganesh Bhandari for his useful comments related to this paper.


\begin{thebibliography}{BLMS}

\bibitem[BK]{BK} S.~Bloch and K.~Kato. $p$-adic \'etale cohomology. \emph{Inst. Hautes \'Etudes Sci. Publ. Math.} No.~63 (1986), 107--152.

\bibitem[BLMS]{BLMS} G.~Bhandari, N.~Lemire, J.~Min\'a\v{c}, and J.~Swallow.  Galois module structure of Milnor $K$-theory in characteristic $p$.  Preprint arXiv:math.NT/0405503v5 (2007).

\bibitem[D]{D} P.~K.~Draxl. \emph{Skew fields}. London Mathematical Society Lecture Note Series 81. Cambridge: Cambridge University Press, 1983.

\bibitem[FK]{FK} I.~Fesenko and M.~ Kurihara, eds. \textit{Invitation to higher local fields (M\"unster, 1999)}. Geometry \& Topology Monographs 3. Somerville, MA: International Press, 2000.

\bibitem[I]{I} O.~Izhboldin.  On $p$-torsion in $K^M_*$ for fields of characteristic $p$.  \textit{Algebraic $K$-theory}. Advances in Soviet Mathematics 4.  Providence, RI: American Mathematical Society, 1991, 129--144.

\bibitem[J]{J} N.~Jacobson. \emph{Finite-dimensional division algebras}. New York: Springer-Verlag, 1996.

\bibitem[K]{K} K.~Kato. Residue homomorphisms in Milnor $K$-theory. \emph{Advanced Studies in Pure Math.} \textbf{2} (1983), 153--172.

\bibitem[Ma]{Ma} B.~Magurn. An algebraic introduction to $K$-theory. Encyclopedia of Mathematics and its Applications 87. Cambridge: Cambridge University Press, 2002.

\bibitem[Mi]{Mi} J.~Milnor. Algebraic $K$-theory and quadratic forms. \textit{Invent.~Math.} {\bf 9} (1970), 318--344.

\bibitem[MSS1]{MSS1} J.~Min\'a\v{c}, A.~Schultz, and J.~Swallow. Galois module structure of $p$th-power classes of cyclic extensions of degree $p^n$. \emph{Proc. London Math. Soc.} {\bf 92} (2006), no.~2, 307--341.

\bibitem[MSS2]{MSS2} J.~Min\'a\v{c}, A.~Schultz, and J.~Swallow. Cyclic algebras and construction of some Galois modules. Preprint arXiv:math.NT/0410536v1 (2004).

\bibitem[P]{P} R.~Pierce. \emph{Associative algebras}.  Graduate Texts in Mathematics 88.  New York: Springer-Verlag, 1982.

\bibitem[SV]{SV} A.~Schofield and M.~Van den Bergh. The index of a Brauer class on a Brauer-Severi variety. \emph{Trans. Amer. Math. Soc.} {\bf 333} (1992), no. 2, 729--739.

\bibitem[T]{T} J. Th\'evenaz. Representations of finite groups in characteristic $p^r$. \emph{J.~Algebra} {\bf 72} (1981), 478--500.

\end{thebibliography}
\end{document}